\documentclass[leqno, 10pt]{amsart}
\usepackage{amsmath}
\usepackage{amssymb, latexsym}
 \newtheorem{proposition}{Proposition}[section]

 \newtheorem*{theorem*}{Theorem}
 
 \newtheorem*{proposition*}{Proposition}

 \theoremstyle{remark}

\newcommand{\op}[1]{\operatorname{#1}}

\newcommand{\acou}[2]{\ensuremath{\langle #1 , #2 \rangle}}

\newcommand{\brak}[1]{\ensuremath{\langle\! #1\!\rangle}}

\newcommand{\C}{\ensuremath{\mathbb{C}}} 
 
\newcommand{\bH}{\ensuremath{\mathbb{H}}} 
\newcommand{\R}{\ensuremath{\mathbb{R}}}

\newcommand{\Rn}{\ensuremath{\R^{2n-1}}}
\newcommand{\URn}{U\times\R^{2n-1}}
\newcommand{\URno}{U\times(\R^{2n-1}\setminus 0)}

\newcommand{\Ca}[1]{\ensuremath{\mathcal{#1}}}

\newcommand{\cE}{\Ca{E}}

\newcommand{\cL}{\ensuremath{\mathcal{L}}}

\newcommand{\cV}{\ensuremath{H}}

\newcommand{\fg}{\mathfrak{g}}

\newcommand{\psivdo}{$\Psi_{H}$DO}
\newcommand{\psivdos}{$\Psi_{H}$DO's}
 
\newcommand{\pvdo}{\ensuremath{\Psi_{\cV}}}

\newcommand{\psido}{$\Psi$DO}

\newcommand{\psidos}{$\Psi$DO's}

\newcommand{\End}{\ensuremath{\op{End}}}

\begin{document}

\title{Comments on: ``Operator $K$-theory for the group ${\rm SU}(n,1)$'' by P.~Julg and G.~Kasparov.}

\author{Rapha\"el Ponge}

\address{Max Planck Institute for Mathematics, Bonn, Germany}

\email{raphaelp@mpim-bonn.mpg.de}

 \begin{abstract}
  In this note we point out and fill a gap in the proof by Julg-Kasparov~\cite{JK:OKTGSU} of the Baum-Connes conjecture with  
 coefficients for discrete subgroups of $\op{SU}(n,1)$. The issue at stake is the proof that the complex powers of the 
 contact Laplacian are element of the Heisenberg calculus. In particular, we explain why we cannot implement  into the setting of the Heisenberg 
 calculus the classical Seeley's approach to complex powers. 
  \end{abstract}
   
 \maketitle  

The aim of this note is to point out  and fill a gap in the proof by Julg-Kasparov~\cite{JK:OKTGSU} of the Baum-Connes conjecture with  
coefficients for discrete subgroups of the complex Lorentz group $G=\op{SU}(n,1)$. The issue is with the proof by Julg and Kasparov that the complex 
powers of the contact Laplacian are pseudodifferential operators in the Heisenberg calculus of Beals-Greiner~\cite{BG:CHM} and Taylor~\cite{Ta:NCMA}. 

To prove this Julg and Kasparov tried to carry out in the Heisenberg setting the classical approach of Seeley~\cite{Se:CPEO} to complex powers of elliptic operators. 
We point a gap in their argument and we further show that we cannot implement Seeley's approach to complex powers into the setting of the Heisenberg calculus. 
Nevertheless, the result about the complex powers of the contact Laplacian can be proved by using the results of~\cite{Po:MAMS1}.  This allows us to 
fill the gap in~\cite{JK:OKTGSU}.

The note is organized as follows. In~Section~\ref{sec:gap} we point out explain a gap in the proof of Julg-Kasparov. In Section~\ref{sec:Heisenberg} we 
give a brief review of the Heisenberg calculus. In Section~\ref{sec:complex-powers} we fill the gap and in Section~\ref{sec:SANCHC} we explain why we actually 
can't implement Seeley's approach in the setting of the Heisenberg calculus. 

\section{A gap in Julg-Kasparov's proof}\label{sec:gap}
The proof by Julg and Kasparov in~\cite{JK:OKTGSU} of the Baum-Connes conjecture with  
coefficients for discrete subgroups of $\op{SU}(n,1)$ can be briefly summarized as follows. 

First, the proof can be reduced to showing that Kasparov's element $\gamma_{G}$ is equal to $1$ 
in the representation ring ${R}(G)$, that is, if $K$ is 
the maximal compact group of $G$ the restriction map ${R}(G)\rightarrow {R}(K)$ is an isomorphism. 

Second, the symmetric space $G/K$ is a complex 
hyperbolic space and under the Siegel map it is biholomorphic to the unit ball $B^{2n}\subset \C^{n}$ and its visual boundary is CR diffeomorphic to 
the unit sphere $S^{2n-1}$ equipped its standard CR structure. Julg and Kasparov further showed that $R(K)$ can be geometrically realized as 
${KK}_{G}(C(\overline{B}^{2n}),\C)$, where $C(\overline{B}^{2n})$ denotes the $C^{*}$-algebra of continuous functions on the closed unit ball 
$\overline{B}^{2n}$ and ${KK}_{G}$ is the equivariant ${KK}$ functor of Kasparov. They then built a Fredholm 
module representing an element $\delta$ in ${KK}_{G}(C^{0}(\overline{B}^{2n}),\C)$ which is mapped to $\gamma$ in $KK_{G}(\C,\C)={R}(G)$ 
under the morphism induced by the map $\overline{B}^{2n}\rightarrow \{\text{pt}\}$. 

The construction of the element $\delta$ in ${KK}_{G}(C(\overline{B}^{2n}),\C)$ involves in a crucial manner  
the contact complex of Rumin~\cite{Ru:FDVC} on the unit sphere $S^{2n-1}$ endowed with its standard 
contact structure. In the contact setting the main geometric operators are not elliptic and the relevant pseudodifferential calculus to deal with them is 
the Heisenberg calculus of Beals-Greiner~\cite{BG:CHM} and Taylor~\cite{Ta:NCMA}. Then for constructing the Fredholm module representing 
$\delta$ Julg and Kasparov had to prove that the complex powers of the contact Laplacian are pseudodifferential operators in the Heisenberg calculus 
(see~\cite[Thm.~5.27]{JK:OKTGSU}). 

In~\cite{Se:CPEO} Seeley settled a general procedure for constructing complex powers of elliptic operators as pseudodifferential operators. Its approach 
relied on constructing asymptotic resolvents in a suitable class of classical \psidos\ calculus with parameter. Accordingly, Julg and Kasparov tried
to construct an asymptotic resolvent for the contact Laplacian 
 in a class of Heisenberg \psidos\ with parameter $\lambda$ in any given angular sector $\Lambda \subset \C \setminus (0, \infty)$ (see~\cite[Thm.~5.25]{JK:OKTGSU}). 
 In particular, in their construction  the symbol in a local chart of the asymptotic resolvent is never defined for $\lambda=0$ and 
 its homogeneous components have meromorphic singularities near $\lambda=0$. 

Now, in order to carry out Seeley's approach for the contact Laplacian we have to be able to integrate the symbol of the asymptotic resolvent with 
respect to the parameter $\lambda$ over a contour $\Gamma$ crossing the value $\lambda=0$. 
This becomes troublesome when in the proof of~\cite[Thm.~5.27]{JK:OKTGSU} 
the authors claim that by their Theorem~5.25 the resolvent of the contact Laplacian belongs to a class of Heisenberg \psidos\ with parameter in a set containing 
$\Gamma$.  In particular, their Theorem~5.25 does not allow them to integrate over 
$\Gamma$ the homogeneous components of symbol of the asymptotic resolvent. This shows that the proof of their Theorem~5.27 is not complete.

%
\section{Heisenberg calculus}\label{sec:Heisenberg}
Let $M^{2n-1}$ be a compact orientable contact manifold with contact hyperplane $H=\ker \theta$, where $\theta$ is a contact form, i.e., 
$d\theta_{|_{H}}$ is 
non-degenerate. The contact condition implies that there is a non-degenerate Levi form 
$\cL: H\times H\rightarrow TM/H$
 such that, for any $x \in M$ and any sections $X'$ and $Y'$ of $H$ we have $\cL_{x}(X'(x),Y'(x))=[X',Y'](x) \ \bmod H_{x}$, i.e., the value of  
 $[X',Y'](x)$ modulo $H_{x}$ only depends on the values at $x$ of $X'$ and $Y'$ and not on their higher order jets. This allows us to define a tangent Lie 
 group bundle $GM$ as the bundle $(TM/H)\oplus H$ equipped with the dilations and group law such that
\begin{gather}
    t.(X_{0}+X')=t^{2}X_{0}+tX', \quad t \in \R,\\
    (X_{0}+X').(Y_{0}+Y')=X_{0}+Y_{0}+\frac{1}{2}\cL(X',Y')+X'+Y',
\end{gather}
for sections $X_{0}$ and $Y_{0}$ of $TM/H$ and sections $X'$ and $Y'$ of $H$. Furthermore, since $H$ is a contact hyperplane $GM$ is in fact a fiber 
bundle of Lie groups with typical fiber the $(2n+1)$-dimensional Heisenberg group $\bH^{2n-1}$.

The Heisenberg calculus of Beals-Greiner~\cite{BG:CHM} and Taylor~\cite{Ta:NCMA} is suitable pseudodifferential calculus to deal with hypoelliptic 
operators on contact manifolds. Its idea, which goes back to Eli Stein, is to construct a class of 
pseudodifferential operators, called \psivdos, whose calculus is modelled on that of convolutions operators on the Heisenberg group. 

Locally \psivdos\ can be described as follows. Let $U \subset \Rn$ be a local chart with a $H$-frame $X_{0},X_{1},\ldots,X_{2n}$ of $TU$, i.e., a 
frame such 
that $X_{1},\ldots,X_{2n}$ span $H$. In the sequel such a chart will be called a Heisenberg chart. In addition, we endow $\Rn$ with the pseudo-norm $ 
\|\xi\|=(\xi_{0}^{2}+\xi_{1}^{4}+\ldots+\xi_{2n}^{4})^{\frac{1}{4}}$ and the dilations $t.\xi=(t^{2}\xi_{0},t\xi_{1},\ldots,t\xi_{2n})$, $t\in \R$.

A Heisenberg symbol of order $m$, $m\in \C$, is a function $p\in C^{\infty}(\URn)$ admitting  an asymptotic expansion $p\sim \sum_{j \geq 
0}p_{m-j}$ with symbols $p_{m-j}\in C^{\infty}(\URno)$ such that $p_{m-j}(x,t.\xi)=t^{m-j}p_{m-j}(x,\xi)$ for any $t>0$. Here the sign~$\sim$ means 
that, for any compact $K \subset U$ and any integer $N$, we have estimates,
\begin{equation}
   |\partial_{x}^{\alpha}\partial_{\xi}^{\beta}(p-\sum_{j<N}p_{m-j})(x,\xi)|\leq C_{KN\alpha\beta}\|\xi\|^{\Re m-N-\brak\beta}, \quad x\in K,\ 
   \|\xi\|\geq 1,
\end{equation}
where we have let $\brak\beta=2\beta_{0}+\beta_{1}+\ldots.+\beta_{2n}$.

Let $\frac{1}{i}X_{j}=\sum_{k} a_{jk}(x)\partial_{x_{k}}$ and set $a(x)=(a_{jk}(x))$. 
A \psivdo\ of order $m$ on $U$ is a continuous operator $P:C^{\infty}_{c}(U)\rightarrow C^{\infty}(U)$ of the form,
\begin{equation}
    Pu(x)=(2\pi)^{-(2n+1)}\int e^{i\acou{x}{\xi}}p(x,a(x)\xi)\hat{u}(\xi)d\xi +Ru(x),
\end{equation}
where $p(x,\xi)$ is a Heisenberg symbol of order $m$ and $R$ is a smoothing operator. 
The class of \psivdos\ of order $m$ is invariant under changes 
of Heisenberg charts and so we can define \psivdos\ of order $m$ on $M$ acting on sections of a vector bundle $\cE$ over $M$. We let $\pvdo(M,\cE)$ denote 
 the class of such operators.

 Let $\fg^{*}M=(TM/H)^{*}\oplus H^{*}$ be the linear dual of the Lie algebra bundle of $GM$. If $P$ is an operator in $\pvdo^{m}(M,\cE)$ then 
 its principal symbol can be invariantly defined as an element $\sigma_{m}(P)(x,\xi)$ of the space $S_{m}(\fg^{*}M,\cE)$ of sections 
$p\in S_{m}(\fg^{*}M\setminus 0,\End \cE)$ such that $p(x,t.\xi)=t^{m}p(x,\xi)$ for any $t>0$ (here $\cE$ is seen as a vector bundle over 
$\fg^{*}M$ using the canonical submersion $\fg^{*}M\rightarrow M$).  

Let $a\in M$ and let $S_{m}(\fg^{*}_{a}M,\cE_{a})$ be the space of functions $p\in C^{\infty}(\fg^{*}_{a}M\setminus 0,\cE_{a})$ which are homogeneous of 
degree $m$. Then under the Fourier transform the convolution of distributions on $G_{a}M$ defines a bilinear product
$*^{a}$ from $S_{m_{1}}(\fg^{*}_{a}M,\cE_{a})\times S_{m_{2}}(\fg^{*}_{a}M,\cE_{a})$ to $S_{m_{1}+m_{2}}(\fg^{*}_{a}M,\cE_{a})$.
This product depends smoothly on $a$, so it gives rise to the product, 
\begin{equation}
    *:S_{m_{1}}(\fg^{*}M,\cE)\times S_{m_{2}}(\fg^{*}M,\cE) \rightarrow S_{m_{1}+m_{2}}(\fg^{*}M,\cE),
    \label{eq:CPCL.product-symbols}
\end{equation}
such that $p_{m_{1}}*p_{m_{2}}(a,\xi)=(p_{m_{1}}(a,.)*^{a}p_{m_{2}}(a,.))(\xi) \ \forall p_{m_{j}}\in S_{m_{j}}(\fg^{*}M,\cE)$.
This product corresponds to the product of \psivdos\ at the level of principal symbols, i.e., we have 
$\sigma_{m_{1}+m_{2}}(P_{1}P_{2})=\sigma_{m_{1}}(P_{1})*\sigma_{m_{2}}(P_{2}) \ \forall P_{j} \in \pvdo^{m_{j}}(M,\cE)$. 

In fact, if in a given local trivializing Heisenberg chart with $H$-frame $X_{0},\ldots,X_{2n-1}$ the operators 
$P\in \pvdo^{m}(M,\cE)$ and $Q\in \pvdo^{m'}(M,\cE)$ have symbols $p\sim \sum p_{m-j}$  and 
$q\sim \sum q_{m'-j}$ then $PQ$ has symbol $r\sim \sum r_{m+m'-j}$, with  
\begin{equation}
     r_{m+m'-j} = \sum_{k+l\leq j} \sum_{\alpha,\beta,\gamma,\delta}^{(j-k-l)}
            h_{\alpha\beta\gamma\delta}  (D_{\xi}^\delta p_{m-k})* (\xi^\gamma 
            \partial_{x}^\alpha \partial_{\xi}^\beta q_{m'-l}), 
            \label{eq:HC.composition}
\end{equation}
where $\underset{\alpha\beta\gamma\delta}{\overset{\scriptstyle{(k)}}{\sum}}$ denotes the sum over all the indices such that 
$|\alpha|+|\beta| \leq \brak\beta -\brak\gamma+\brak\delta = k$ and $|\beta|=|\gamma|$, and the functions 
$h_{\alpha\beta\gamma\delta}(x)$'s are  polynomials in the derivatives of the coefficients of 
the vector fields $X_{0},\ldots,X_{2n-1}$. It follows from this that we can construct a parametrix of $P$ in $\pvdo^{-m}(M,\cE)$ if, and only if, its 
principal symbol $\sigma_{m}(P)$ is invertible with respect to the product $*$. 

As $G_{a}M$ is not Abelian the product $*^{a}$ is not anymore the pointwise product of functions. Therefore, if $p_{m-j}\in S_{m_{j}}(\fg^{*}_{a}M)$, 
$j=1,2$, 
then the computation of  $p_{m_{1}}*p_{m_{2}}$ at a point $\xi\in \fg^{*}M\setminus 0$
requires the knowledge of the values  $p_{m_{1}}$ and $p_{m_{2}}$ at all the points $\xi'$ of $\fg^{*}_{a}M$. It follows that the 
product~(\ref{eq:CPCL.product-symbols})  
for Heisenberg symbols is not microlocal, i.e., it cannot be localized with respect to the $\xi$-variable. 

As we will explain in Section~\ref{sec:SANCHC} this lack of microlocality prevents us  from carrying out Seeley's approach to complex powers. 
Nevertheless, complex powers of hypoelliptic \psidos\ were dealt with in~\cite{Po:MAMS1} by relying on a pseudodifferential representation of the heat kernel, 
instead of using a pseudodifferential representation of the 
resolvent as in Seeley's approach.  In particular, assuming $M$ and $\cE$ endowed with a Riemannian metric and a Hermitian metric, we have:

%

\begin{proposition}[\cite{Po:MAMS1}]\label{prop:Heisenberg.complex-powers}
    Let $P:C^{\infty}(M,\cE)\rightarrow C^{\infty}(M,\cE)$ be a differential operator~$\geq 0$ of Heisenberg order $m$ such that $\sigma_{m}(P)$ is 
    invertible. Then, for any $s \in \C$, the power $P^{s}$ defined by $L^{2}$-functional calculus is a \psivdo\ of order $ms$.
\end{proposition}

\section{Complex powers of the contact Laplacian}\label{sec:complex-powers}
 Let $M^{2n-1}$ be a compact orientable contact manifold with contact hyperplane $H=\ker \theta$, where $\theta$ is a contact form and 
 let $J$ be an almost complex structure on $H$ such that we have $d\theta(X,JX)>0$ for 
any section $X$ of $H\setminus 0$. We then can endow $M$ with the Riemannian metric $g_{\theta,J}=\theta^{2}+d\theta(.,J.)$. In addition,
we let $X_{0}$ be the Reeb vector 
field associated to $\theta$, so that $\theta(X_{0})=1$ and $\iota_{X_{0}}d\theta=0$. 

The splitting $TM=\R X_{0}\oplus H$ allows us to identify  
$H^{*}$ with the annihilator of $X_{0}$ in $T^{*}M$ and to identify $\Lambda^{k}_{\C}H^{*}$ with $\ker \iota_{X_{0}}$, so that we get the 
orthogonal splitting $ \Lambda^{*}_{\C}TM=(\bigoplus_{k=0}^{2n} \theta\wedge 
\Lambda^{k}_{\C}H^{*}) \oplus (\bigoplus_{k=0}^{2n}\Lambda^{k}_{\C}H^{*})$. 
 If $\eta\in C^{\infty}(M,\Lambda^{k}_{\C}H^{*})$ then we have
$d\eta= \theta \wedge \cL_{X_{0}}\eta+d_{b}\eta$,
where $d_{b}\eta$ is the component of $d\eta$ in $\Lambda^{k}_{\C}H^{*}$. This does not provide us with a complex, for we have 
$d_{b}^{2}=-\cL_{X_{0}}\varepsilon(d\theta)$, where $\varepsilon(d\theta)$ denotes the exterior multiplication 
by $d\theta$. 

The contact complex of Rumin~\cite{Ru:FDVC} arises as a modification $d_{b}$ and $\Lambda^{k}_{\C}H^{*}$ to 
get a complex of horizontal differential forms. Let $\Lambda^{*}_{1}:=\ker \iota(d\theta)\cap \Lambda^{*}_{\C}H^{*}$ and 
$\Lambda^{*}_{2}:=\ker \varepsilon(d\theta) \cap \Lambda^{*}_{\C}H^{*}$. Then $C^{\infty}(M,\Lambda^{*}_{2})$
is closed under $d_{b}$ 
and annihilates $d_{b}^{2}$ and $C^{\infty}(M,\Lambda^{*}_{2})$ is closed under $d_{b}^{*}$ and annihilated by $(d_{b}^{*})^{2}$, so that we get 
two complexes. However, 
since $d\theta$ is nondegenerate on $H$ the operator $\varepsilon(d\theta):\Lambda^{k}_{\C}H^{*}\rightarrow \Lambda^{k+2}_{\C}H^{*}$  is 
injective for $k\leq n-2$ and surjective for $k\geq n$ and so $\Lambda_{2}^{k}=\{0\}$ for $k\leq n-1$ and $\Lambda_{1}^{k}=\{0\}$ for $k\geq n$. 
Therefore, we only have two halves of complexes. 

As observed by Rumin~\cite{Ru:FDVC} we get a full complex by connecting the two halves by means of the 
differential operator $D_{R;n-1}:=\cL_{X_{0}}+d_{b;n-2}\varepsilon(d\theta)^{-1}d_{b;n-1}$ acting on $C^{\infty}(M,\Lambda_{\C}^{n-1}H^{*})$,
where $\varepsilon(d\theta)^{-1}$ is the inverse of $\varepsilon(d\theta):\Lambda^{n-2}_{\C}H^{*}\rightarrow \Lambda^{n}_{\C}H^{*}$. Therefore, if 
we let $\pi_{1}\in C^{\infty}(M,\Lambda_{\C}^{*}H^{*})$ be the orthogonal projection onto $\Lambda_{1}$ then we 
have the complex, 
 \begin{equation}
     C^{\infty}(M)\stackrel{d_{R;0}}{\rightarrow}
     \ldots
     C^{\infty}(M,\Lambda^{n-2})\stackrel{D_{R;n-1}}{\rightarrow} C^{\infty}(M,\Lambda^{n-1}) 
     \ldots \stackrel{d_{R;2n-3}}{\rightarrow} C^{\infty}(M,\Lambda^{2n-2}).
      \label{eq:contact-complex}
 \end{equation}
where $d_{R;k}$ agrees with $\pi_{1}\circ d_{b}$ for $k=0,\ldots,n-2$ and  with $d_{R;k}=d_{b}$ otherwise. 

The contact Laplacian is defined as follows. In degree $k\neq n$ this is the differential operator 
$\Delta_{R;k}:C^{\infty}(M,\Lambda^{k})\rightarrow C^{\infty}(M,\Lambda^{k})$ such that
\begin{equation}
    \Delta_{R;k}=\left\{
    \begin{array}{ll}
        (n-1-k)d_{R;k-1}d^{*}_{R;k}+(n-k) d^{*}_{R;k+1}d_{R;k},& \text{$k=0,\ldots,n-2$},\\
         (k-n)d_{R;k-1}d^{*}_{R;k}+(k-n+1) d^{*}_{R;k+1}d_{R;k},& \text{$k=n,\ldots,2n$}.
         \label{eq:contact-Laplacian1}
    \end{array}\right.
\end{equation}
For $k=n-1$ we have the differential operators $\Delta_{R;n-1,j}:C^{\infty}(M,\Lambda_{j}^{n})\rightarrow C^{\infty}(M,\Lambda^{n}_{j})$, $j=1,2$, 
given by the formulas, 
\begin{gather}
    \Delta_{R;n-1,1}= (d_{R;n-2}d^{*}_{R;n-1})^{2}+D_{R;n-1}^{*}D_{R;n-1}, \\   \Delta_{R;n-1,2}=D_{R;n-1}D_{R;n-1}^{*}+  (d^{*}_{R;n}d_{R;n-1}).
\end{gather}
Notice that $\Delta_{R;k}$, $k \neq n-1$, is a differential operator of (Heisenberg) order $2$, while $\Delta_{R;n-1,j}$, $j=1,2$, is a differential operator 
of (Heisenberg) order $4$. 

It has been shown by Rumin~\cite{Ru:FDVC} that  the contact Laplacian is hypoelliptic. In fact, we have:

\begin{proposition}[\cite{JK:OKTGSU}, \cite{Po:MAMS1}]
    The principal symbols $\sigma_{2}(\Delta_{R;k})$, $k \neq n-1$, and $\sigma_{4}(\Delta_{R;n-1,j})$, $j=1,2$, are invertible with the respect to the 
    product~(\ref{eq:CPCL.product-symbols}) for Heisenberg symbols.
\end{proposition}

Combining this with Proposition~\ref{prop:Heisenberg.complex-powers} then gives:

\begin{proposition}[\cite{Po:MAMS1}]
    For any $s\in \C$ the powers $\Delta_{R;k}^{s}$, $k\neq n-1$, and $\Delta_{R;n-1,j}^{s}$, $j=1,2$, defined by $L^{2}$-functional calculus are 
    \psivdos\ of degree $2s$ and $4s$ respectively.
\end{proposition}
This fills the gap in~\cite{JK:OKTGSU} alluded to in Section~\ref{sec:gap}.

\section{Heisenberg calculus and Seeley's approach}
\label{sec:SANCHC}
In this last section we  explain why the lack of microlocality of the Heisenberg calculus actually prevents us from implementing into 
this setting Seeley's approach to complex powers. 

\subsection{Seeley's approach to complex powers}
 Let us briefly recall the approach of Seeley~\cite{Se:CPEO} to complex powers (see also~\cite{GS:WPPOAPSBP}, \cite{Sh:POST}). 
 To simplify the exposition we let $M^{n}$ be a compact Riemannian manifold 
 equipped and let $\Delta:C^{\infty}(M)\rightarrow C^{\infty}(M)$ be a second order positive elliptic differential operator with principal 
 symbol  $p_{2}(x,\xi)>0$.  Then for $\Re s<0$ we have: 
 \begin{gather}
    \Delta^{s}=\frac{i}{2\pi} \int_{\Gamma_{r}} \lambda^{s}(\Delta-\lambda)^{-1}d\lambda, 
     \label{eq:SA.complex-powers-definition}\\
\Gamma_{r}=\{ \rho e^{i\pi}; \infty <\rho\leq r\}\cup\{ r e^{it}; 
\theta\geq t\geq \theta-2\pi \}\cup\{ \rho e^{-i\pi};  r\leq \rho\leq \infty\},
\label{eq:SA.contour}
 \end{gather}
 where $r>0$ is small enough so that non nonzero eigenvalue of $\Delta$ lies in $(0,r]$. 
 
 To show that the formula above defines a \psido\ Seeley constructs an asymptotic resolvent $Q(\lambda)$ as a parametrix for $\Delta-\lambda$ in a 
 suitable 
 \psido\ calculus 
 with parameter. More precisely, let $\Lambda \subset \C\setminus 0$ be an open angular sector $\theta <\arg \lambda<\theta'$ with $0<\theta<\pi<\theta'<2\pi$. 
 In the sequel we will say that a subset $\Theta \subset 
 [\R^{n}\times \C]\setminus 0$ is \emph{conic} when for any $t>0$ and any $ (\xi,\lambda)\in \Theta$ we have $(t\xi,t^{2}\lambda)\in \Theta$.
 For instance the subset $\R^{n}\times \Lambda \subset [\R^{n}\times \C]\setminus 0$ is conic. 
 
 Let $U \subset \R^{n}$ be a local chart for $M$. Then in $U$ the asymptotic resolvent has a symbol of the form
$q(x,\xi;\lambda) \sim \sum_{j\geq 0}q_{-2-j}(x,\xi;\lambda)$,
 where $\sim$ is taken in a suitable sense (see~\cite{Se:CPEO}) and there exists an open conic subset $\Theta\subset [\R^{n}\times \C]\setminus 0$ containing 
 $\R^{n}\times \Lambda$ such that each symbol $q_{-2-j}(x,\xi;\lambda)$ is smooth on $U\times \Theta$ and satisfies
 \begin{equation}
     q_{-2-j}(x,t\xi;t^{2}\lambda)=t^{-2-j}q(x,\xi;\lambda) \quad \forall t>0.
      \label{eq:SA.homogeneity}
 \end{equation}
 
 If $p(x,\xi)=\sum_{j=0}^{2}p_{2-j}(x,\xi)$ denotes the symbol of $\Delta$ in the local chart $U$ then $ q(x,\xi;\lambda)$ is such that
 $1\sim  (p(x,\xi)-\lambda)q(x,\xi;\lambda)+\sum_{\alpha \neq 0}\frac{1}{\alpha !}\partial_{\xi}^{\alpha}p(x,\xi)D_{x}^{\alpha}q(x,\xi;\lambda)$,
 from which we get
 \begin{gather}
     q_{-2}(x,\xi;\lambda)= (p_{2}(x,\xi)-\lambda)^{-1},
     \label{eq:SA.symbol-resolvent1}\\
      q_{-2-j}(x,\xi;\lambda)=-q_{-2}(x,\xi;\lambda)\!\! \sum_{\substack{|\alpha|+k+l=j,\\ l\neq j}} \!\! 
      \frac{1}{\alpha !}\partial_{\xi}^{\alpha}p_{2-k}(x,\xi)D_{x}^{\alpha}q_{-2-l}(x,\xi;\lambda).
     \label{eq:SA.symbol-resolvent2}
 \end{gather}
 
Set $\rho=\inf_{x \in U}\inf_{|\xi|=1}p_{2}(x,\xi)$.  Possibly by shrinking $U$ we may assume $\rho>0$. Let 
$    \Theta =[\R^{n}\times \Lambda] \cup\{(\xi;\lambda)\in \R^{n} \times \C; 0\leq |\lambda|<\rho|\xi|^{2}\}$.
 Then the formulas~(\ref{eq:SA.symbol-resolvent1}) and (\ref{eq:SA.symbol-resolvent2}) show that each symbol $q_{-2-j}(x,\xi;\lambda)$ is well defined and 
 smooth on $U\times \Theta$ and is homogeneous in 
 the sense of~(\ref{eq:SA.homogeneity}). Furthermore, it is analytic with respect to $\lambda$. Therefore, for $\Re s<0$ we define a smooth function on $U\times 
( \R^{n}\setminus 0)$ by letting 
\begin{equation}
    p_{s,ms-j}(x,\xi)= \frac{i}{2\pi} \int_{\Gamma_{\xi}} \lambda^{s}q_{-2-j}(x,\xi;\lambda)d\lambda,
\end{equation}
where $\Gamma_{\xi}$ is the contour $\Gamma_{r}$ in~(\ref{eq:SA.contour}) with $r=\frac{1}{2}\rho|\xi|^{2}$. Moreover, one can check that 
$p_{s,ms-j}(x,t\xi)=t^{2s-j}p_{s,ms-j}(x,\xi)$ for any $t>0$, i.e., $p_{s,ms-j}(x,\xi)$ is a homogeneous symbol of degree $ms-j$. 

It can also be shown that on the chart $U$ the operators $\Delta^{s}$ is a \psido\ with symbol $p_{s} \sim \sum_{j \geq 0}p_{s,ms-j}(x,\xi)$.  This is true on any 
local chart  and we can check 
that the Schwartz kernel of $\Delta^{s}$ is smooth off the diagonal of $M\times M$, so we see that $\Delta^{s}$ is a \psido\ of order $ms$ for $\Re s<0$. 

Finally, for $\Re s\geq 0$ and $k$ integer~$>\Re s$ we have $\Delta^{s}=\Delta^{k}\Delta^{s-k}$, so since 
$\Delta^{k}$ is a differential operator of order $k$ and  
 $\Delta^{s-k}$ is a \psido\ of order $m(s-k)$ we see that $\Delta^{s}$ is a \psido\ of order $ms$. Hence $\Delta^{s}$ is a \psido\ of order $ms$ for any $s \in \C$. 

\subsection{Obstruction to Seeley's approach} 
Let us now explain why we cannot carry out Seeley's approach within the framework of the Heisenberg calculus. We will explain this in the special 
case of the contact Laplacian $\Delta_{R;0}$ acting on the functions on a compact orientable contact manifold $(M^{2n-1},H)$ as in 
Section~\ref{sec:complex-powers}.

 In order to carry out Seeley's approach for $\Delta_{R;0}$ we have to construct an asymptotic resolvent  in a class of \psivdos\ with 
parameter associated to an angular sector $\Lambda \subset \C\setminus [0,\infty)$ as above and given in a local Heisenberg chart 
$U \subset \Rn$ by parametric symbols, 
$q(x,\xi;\lambda) \sim \sum_{j\geq 0}q_{-2-j}(x,\xi;\lambda)$, 
where $\sim$ is taken in a suitable sense and there exists an open conic subset $\Theta\subset [\R^{n}\times \C]\setminus 0$ containing 
 $\R^{2n-1}\times \Lambda$ such that each symbol $q_{-2-j}(x,\xi;\lambda)$ is smooth on $U\times \Theta$ and satisfies 
 $q_{-2-j}(x,t.\xi;t^{2}\lambda)=t^{-2-j}q(x,\xi;\lambda)$ for any $t>0$. 
 
 If we let $p(x,\xi)=\sum p_{2-j}(x,\xi)$ be the symbol of $\Delta_{R;0}$ in the Heisenberg chart, 
 then by~(\ref{eq:HC.composition}) we have
  \begin{equation}
     1 \sim \sum_{j \geq 0}  \sum_{k+l\leq j} \sum_{\alpha,\beta,\gamma,\delta}^{(j-k-l)}
             h_{\alpha\beta\gamma\delta}(x)  (D_{\xi}^\delta p_{2-k})* (\xi^\gamma 
             \partial_{x}^\alpha \partial_{\xi}^\beta q_{-2-l})(x,\xi;\lambda),
  \end{equation}
 from which we get
\begin{gather}
    q_{-2}(x,\xi;\lambda)=(p_{2}-\lambda)^{*-1}(x,\xi;\lambda)
    \label{eq:Gap.symbol-resolvent1}  \\
    q_{-2-j}(x,\xi;\lambda)=-\sum_{\substack{k+l\leq j,\\ l\neq j}} \sum_{\alpha,\beta,\gamma,\delta}^{(j-k-l)}
            h_{\alpha\beta\gamma\delta}(x)  q_{-2}*(D_{\xi}^\delta p_{2-k})* (\xi^\gamma 
            \partial_{x}^\alpha \partial_{\xi}^\beta q_{-2-l})(x,\xi;\lambda),
    \label{eq:Gap.symbol-resolvent2} 
\end{gather}
where $(p_{2}-\lambda)^{*-1}$ denotes the inverse of $p_{2}-\lambda$ with respect to the product $*$. 

If $q_{1}(x,\xi;\lambda)$ and $q_{2}(x,\xi;\lambda)$ are two homogeneous Heisenberg symbols with parameter then the product $q_{1}$ and $q_{2}$ 
should be defined as 
\begin{equation}
    q_{1}*q_{2}(x,\xi;\lambda)=[q_{1}(x,.;\lambda)*^{x}q_{2}(x,.;\lambda)](\xi).
\end{equation}
As mentioned in Section~\ref{sec:Heisenberg} the definition of $[q_{1}(x,.;\lambda)*^{x}q_{2}(x,.;\lambda)](\xi)$ depends on all the values of 
$q_{1}(x,\xi';\lambda)$ and $q_{2}(x,\xi';\lambda)$ as $\xi'$ ranges over $\Rn\setminus 0$. For a parameter 
$\lambda>0$ the symbols $q_{1}(x,\xi;\lambda)$ and 
$q_{2}(x,\xi;\lambda)$ are only defined for $\xi$ in $\{\xi; (x,\xi;\lambda)\in \Theta\}$ which does not agree with $\Rn\setminus 0$, so we cannot 
define $q_{1}*q_{2}(x,\xi;\lambda)$ for $\lambda>0$. Therefore, the formula~(\ref{eq:Gap.symbol-resolvent2}) does not make sense for $\lambda>0$. 

All this shows that the non-microlocality of the Heisenberg calculus prevents us from implementing Seeley's approach into the setting of the Heisenberg 
calculus. As previously mentioned the results of~\cite{Po:MAMS1} allows us to deal with complex powers in case of positive differential operators with invertible 
principal symbols. It also possible to construct complex powers for more general hypoelliptic \psidos\ in the spirit of Seeley's approach by replacing the use 
of a homogeneous  symbols with parameter by  \emph{almost} homogeneous symbols with parameter (see~\cite{Po:PhD}, \cite{Po:CPDE1}).  

%
%

\end{document}